\documentclass[11pt]{amsart}
\usepackage{amssymb}

\usepackage{graphicx}
\usepackage{amsmath}

\newtheorem{theorem}{Theorem}

\newtheorem{corollary}[theorem]{Corollary}

\newtheorem{definition}[theorem]{Definition}
\newtheorem{example}[theorem]{Example}

\newtheorem{lemma}[theorem]{Lemma}

\newtheorem{proposition}[theorem]{Proposition}
\newtheorem{remark}[theorem]{Remark}

\begin{document}
\title[BCOV]{Generalized Hodge Metrics and BCOV
torsion on Calabi-Yau Moduli}
\date{June 1, 2003}
\begin{abstract}
We establish an unexpected relation among
the Weil-Petersson metric, the generalized Hodge metrics and
the BCOV torsion. Using this relation,
we prove that certain kind of moduli spaces
of polarized Calabi-Yau manifolds do not admit
complete subvarieties. That is, there is no complete
family for certain
class of polarized Calabi-Yau manifolds. We
also give an estimate
of the complex Hessian of the BCOV torsion using the relation.
After establishing a degenerate version of the Schwarz Lemma
of Yau, we prove that  the complex Hessian of
the BCOV torsion is bounded by the Poincar\'e metric.
\end{abstract}

\author{Hao Fang} \author{Zhiqin Lu}
 \subjclass{Primary: 53A30;
Secondary: 32C16} \keywords{Calabi-Yau Manifold, analytic torsion,
 moduli space} \email[Hao Fang]{haofang@cims.nyu.edu}
\email[Zhiqin Lu]{zlu@math.uci.edu}
 \thanks{The first author is partially supported by a grant from the
New York University Research Challenge Fund Program;
the second author is partially supported by the NSF grant DMS
0204667 and the Alfred P. Sloan Research Fellowship.}
\maketitle \tableofcontents \pagestyle{myheadings}
\markboth{Generalized Hodge metrics and BCOV torsion}{H. Fang
and Z. Lu}

\newcommand{\ka}{K\"ahler }
\newcommand{\ii}{\sqrt{-1}}
\newcommand{\eps}{\varepsilon}

\section{Introduction}

A Calabi-Yau manifold is a smooth K\"ahler manifold with trivial
canonical bundle and fundamental group. Moduli space of polarized Calabi-Yau manifolds
(Calabi-Yau moduli) is the object to
study in Mirror Symmetry, hence the focal point of intensive
studies in areas of mathematical physics, algebraic geometry,
differential geometry and number theory. For the general reference
of  Mirror Symmetry and related topics, see the book of Cox
and Katz~\cite{CoxKatz} and the
recent survey paper of Todorov~\cite{To1}.

In this paper we study the differential geometry of Calabi-Yau
moduli. The starting point is  the celebrated theorem of Yau,
which establishes the existence of the Ricci flat K\"ahler metric
for a fixed polarization. Differential geometric objects of
Calabi-Yau moduli are usually  constructed from global
algebro-geometric or analytic properties of the Calabi-Yau
manifolds, for example, variation of Hodge structure, or spectral
properties of the Ricci-flat metrics. In this paper, we focus on
the Weil-Petersson metric, the generalized Hodge metrics, the BCOV
torsion and their relations.

From differential geometric point of view, Calabi-Yau moduli  is
amicable. It is smooth (or at worst with quotient singularities,
so that it is a smooth Deligne-Mumford stack) and its local
uniformization is an integral submanifold of the
horizontal distribution of the classifying space of the variation
of Hodge structure. The curvature of the
first Hodge bundle is positive definite. Hence it defines the
so-called Weil-Petersson metric. This metric initiates the study of the
moduli space in terms of differential geometry. For computing the curvatures
of the Weil-Petersson metrics (in  more general settings) and its applications, see the works of Siu~\cite{Siu}
and Schumacher~\cite{Sch}.

In~\cite{Lu3,Lu5}, the second author introduced a new
metric, the Hodge metric, on the Calabi-Yau moduli, mainly
inspired from the theory of variation of Hodge structure. Both the
Hodge metric and the Weil-Petersson metric are K\"ahler orbifold
metrics. However, the Hodge metric enjoys better curvature
properties:  it has non-positive bisectional curvatures, and
furthermore, its
holomorphic sectional curvature and  Ricci curvature
are negative
and bounded away from zero. This  clearly is  not the case for
the Weil-Petersson metric, as shown in the example of
Calabi-Yau quintics~\cite[page 65]{COGP}. Thus, in terms of
differential geometry, the Hodge metric is better than  the
Weil-Petersson metric on Calabi-Yau moduli.

As natural K\"ahler metrics on a given Calabi-Yau moduli,
the
Weil-Petersson metric and the Hodge metric are closely
related.  Both can be realized as curvature forms of various
combinations of the Hodge bundles, in the sense of
Griffiths~\cite{Gr}. There are explicit relations between the
two metrics for K3 surfaces, Calabi-Yau three and four-folds
(see Theorem~\ref{thm22}).  In higher dimensions,  the
concept of the Weil-Petersson geometry was introduced
in~\cite{LS-1}.

However, in defining the Hodge metric, only those
Hodge bundles of the middle
dimensional primitive cohomology groups are considered.
It turns out
that considering the whole Hodge bundles would be more natural.
This is indeed the case.
The universal deformation space of Calabi-Yau manifolds can
also be viewed as horizontal slices of the classifying space
of any degree primitive cohomology groups,
even though  the period maps fail to be immersive in general.
In this paper, we introduce pseudo-metrics on  the classifying
space for $ H^{k}$ where $k$ may  not be the (complex)
dimension of the Calabi-Yau manifolds. These metrics are
called generalized Hodge metrics. The positive definiteness of
the generalized Hodge metrics is lost due to the possible degeneracy of
the corresponding horizontal slices. Nevertheless,
good ``curvature" properties of the Hodge metric still hold for
these generalized Hodge metrics. See Appendix A for the
precise statements.

We now turn to the other geometric object that
will be studied in
this paper: the BCOV torsion. First introduced by and named
after Bershadsky-Cecotti-Ooguri-Vafa~\cite{BCOV2,BCOV1}, the BCOV
torsion is a smooth function on the Calabi-Yau moduli. It
is defined as:
\begin{equation} T=\prod_{1\leq p,q\leq n}(\det \Delta
_{p,q}^{\prime })^{(-1)^{p+q}pq},
\end{equation}
where $\Delta _{p,q}$ is the $\overline{\partial }$-Laplace operator
on $(p,q)$ forms
with
respect to the Ricci-flat metric on a fiber;
$\Delta _{p,q}^{\prime }$ represents the non-singular part of
$\Delta_{p,q}$; the determinant
is taken in the sense of zeta function regularization.

In physics literature, $T$ was
first introduced as the stringy genus one partition function of
$N=2$ SCFT. It was computed using Physics Mirror Symmetry and
was used to predict the number of holomorphic elliptic curves
embedded in certain Calabi-Yau manifolds.

Due to its central role in the Physics Mirror Symmetry, we are
interested in the analysis of  the BCOV torsion mathematically.
$T$ is  a spectral invariant of the Ricci-flat metrics.
The Weil-Petersson metric and the generalized Hodge metrics,
which are defined by using the variation of Hodge
structure, assume no apparent links to the function $T$. The
main result of this paper is the following surprising relation
between the Weil-Petersson metric and
the generalized Hodge metrics by the
BCOV torsion.

\begin{theorem}\label{main}
Let $\omega_{WP}$, $\omega_H$ and $\omega_{H^i}$ be the K\"ahler form of the
Weil-Petersson metric, the Hodge metric and the generalized Hodge metrics,
respectively (See \S
\ref{section2} for the definition). Then
\begin{equation}\label{markertest}
\sum\limits_{i=1}^{n}(-1)^{i}\omega _{H^{i}}
-\frac{\sqrt{-1}}{2\pi}\partial\bar\partial\log
T=\frac{\chi_{Z}}{12}\omega _{WP},
\end{equation}
where $\chi_Z$  is the Euler characteristic number of $Z$.
In particular, if the
Calabi-Yau manifold is primitive
(See \S 3 for the definition), then
\begin{equation}\label{marker2}
\omega_H=(-1)^n(\frac{\sqrt{-1}}{2\pi}\partial\bar\partial\log
T+
\frac{\chi_{Z}}{12}\omega _{WP})).
\end{equation}
\end{theorem}

 As the first
application, we have the following

\begin{corollary}
\label{mo}If $%
N\subset \mathcal{M}$ is a $k$-dimensional complete
subvariety
of $\mathcal M$ where $\mathcal M$ is  the moduli space of a
primitive Calabi-Yau manifold, then
\begin{equation}
\text{Vol}_{H}(N)=[\frac{(-1)^{n}}{12}\chi
_{Z}]^{k}\text{Vol}_{{\normalsize  WP}}(N),
\end{equation}
where $\text{Vol}_{H}(N)$ and $\text{Vol}_{{\normalsize
WP}}(N)$ are the volumes of $N$ with respect to the Hodge and
the Weil-Petersson metrics.
\end{corollary}

In Corollary \ref{mo}, the BCOV torsion does not
appear explicitly.\ Even in dimension 3 and 4, where the Hodge
metric can be expressed explicitly by the Weil-Petersson metric
and its Ricci curvature (\cite{Lu5},~\cite{LS-1}), this volume
identity is new. One of the notable consequences of
the volume identity is the following

\begin{corollary}\label{mi}
Assume that a polarized  Calabi-Yau manifold $Z$ is
primitive, and that
$(-1)^{n+1}\chi _{Z}>-24$.
Let $\mathcal M$ be the moduli space of $Z$.
Then there exists no complete curve in
$\mathcal{M} $; hence, there exists no projective subvariety
of
$\mathcal{M}$ (of positive dimensions). In particular,
$\mathcal M$ is not compact.
\end{corollary}

This corollary is  purely algebro-geometric.
Primitive Calabi-Yau manifolds include
interesting examples like  Calabi-Yau three-folds and
Calabi-Yau hyper-surfaces in projective spaces. It would be
interesting to see a direct proof of the result
without using differential geometry.

The second application of Theorem~\ref{main} is
on the
asymptotic behavior of the complex Hessian of the BCOV
torsion.

\begin{corollary}\label{cor14}
Let $\Delta$ and $\Delta^*$ be the unit disk and the punctured
unit disk of  $\mathbb C$ respectively. Let $(\Delta^*)^l
\times\Delta^{m-l}$ be the parameter space of a family of
Calabi-Yau manifolds. Then the BOCV torsion $T$, which is a
smooth function on $(\Delta^*)^l
\times\Delta^{m-l}$, satisfies
\[
-C\omega_P<\frac{\sqrt{-1}}{2\pi}\partial\bar\partial
\log T< C\omega_P,
\]
where $C$ is a constant and
 $\omega_P$ is the Poincar\'e metric, defined as
\begin{equation}\label{eq15}
\omega_P=\sum_{i=1}^l\ii\frac{1}{
|z_i|^2(\log\frac{1}{|z_i|})^2} dz_i\wedge d\bar z_i
+\sum_{i=l+1}^m
\ii dz_i\wedge d\bar z_i.
\end{equation}
\end{corollary}

Using Theorem~\ref{main}, the proof of the above corollary is
reduced to the fact that the generalized Hodge
metrics are bounded by the Poincar\'e metric. In the
Appendix, we establish Theorem~\ref{A1}, a degenerate version
of the Schwarz Lemma of Yau~\cite{Y3}, of which
Corollary~\ref{cor14} is a direct consequence.

The relation among the Weil-Petersson metric, the generalized
Hodge metric, and the BCOV torsion on  Calabi-Yau moduli
presented in this paper is quite delicate and our
understanding of it is far from being complete.\ Because
of the
physics background of the BCOV torsion, it is very  likely
that some deeper relations can be used to explain the current
coincidence. It is also expected that these constructions
will produce new modular forms on various moduli spaces, as
 the previous works of Yoshikawa~\cite{Yoshi-3, Yoshi-1}
indicated.

We shall  proceed to study  the asymptotic
behavior of the BCOV torsion near the boundary of
Calabi-Yau moduli,
and  the BCOV prediction of  counting the  rational
curves. The results will be the subject of an upcoming paper.

\bigskip

{\bf Acknowledgement.} Both authors are indebted
to G.\ Tian for
introducing the papers~\cite{BCOV2,BCOV1} to them and being
interested in this work.   They also thank J.-P.\ Bismut,
S.-Y. A. Chang, P. Sarnak and P.\ Yang for their interest in this work.
The first author also thanks A.\ Ching for discussion.

\section{Generalized Hodge Metrics}\label{section2}

\setcounter{equation}{0} \setcounter{theorem}{0}

Let $Z$ be a Calabi-Yau manifold and let $l$ be an ample line
bundle over $Z$. The pair $(Z,l)$ is called a polarized
Calabi-Yau manifold. The (coarse) moduli space $\mathcal M$
exists and is constructed as follows: first, choose a large
integer $k$ such that $l^k$ is very ample. In this way $Z$ is
embedded into a complex projective space $\mathbb C P^N$. Let
 $\frak{Hilb}(Z)$  be the Hilbert scheme of $Z$, which  is
a compact complex variety. The group $G=PSL(N+1,\mathbb C)$
acts on
$\frak{Hilb}(Z)$ and the moduli space $\mathcal M$ is the quotient
of the stable points of  $\frak{Hilb}(Z)$ by the group $G$.

By the smoothness theorem of Tian~\cite{T1} (see also
Todorov~\cite{To}), the deformation of the complex structures of
Calabi-Yau manifold is unobstructed. That is, the universal
deformation space (Kuranishi space) is smooth. On the other hand,
due to the existence of finite automorphism, the moduli space for
polarized Calabi-Yau manifolds may have quotient singularities.
Thus in general, the moduli space $\mathcal M$ is a complex
orbifold, or a smooth Deligne-Mumford stack.

For the local geometry of the
moduli space, the above possible singularities may
never be a problem. We can always pass through
a finite covering and assume that locally, the moduli
space is a smooth complex manifold.

Let $Z$ be a generic polarized Calabi-Yau manifold. There exists
the universal family $\frak X$ such that it is parametrized by the
moduli space $\mathcal M$:

\begin{equation}
\begin{array}{ccc}
Z & \overset{i}{\rightarrow } & \frak{X} \\
&  & \,\,\,\,\,\downarrow {\small \pi } \\
&  & \mathcal{M}
\end{array}
,
\end{equation}

Assume that $U$ is an open neighborhood of a point $t\in\mathcal
M$. By the Kodaira-Spencer deformation theory, there is an
isomorphism
\begin{equation}
\iota:
T_{t}U\cong H^{1}(Z_{t},\Theta_t),
\end{equation}
where $Z_t$ is the fiber of $\frak X\overset{\pi}
{\rightarrow}\mathcal M$ at
$t$, and $\Theta_t$ is the holomorphic tangent
bundle of $Z_t$.\footnote{
If $t$ is a singular point of $\mathcal M$,
then we should replace $U$ by $\hat U$, the local
uniformization of $U$.}

\newcommand{\pa}{\partial}

Let $(t_1,\cdots, t_m)$ be a local holomorphic
coordinate system of $\mathcal M$. Then
$\iota(\frac{\pa}{\pa t_i} )\in H^1(Z_t,\Theta_t)$.
We define a Hermitian inner product on $T_t\mathcal
M$ by
\[
\left(\frac{\pa}{\pa t_i},\frac{\pa}{\pa \bar t_j}
\right)_{WP}
=\int_{Z_t} A_{i\bar\beta}^\alpha\cdot
\overline{A_{j\bar\delta}^\gamma} g^{\delta\bar\beta}
g_{\alpha\bar\gamma} dV_{Z_t},
\]
where $A_i=A_{i\bar\beta}^\alpha\frac{\partial}{\partial
t_i}\otimes d\bar t^\beta$, $(i=1,\cdots,m)$ are the harmonic
representation of
$\iota(\frac{\partial}{\partial t_i})$. This inner product on
each $T_t U$ for $t\in\mathcal  M$ gives
a Hermitian metric on the moduli space $\mathcal M$, which is
called the Weil-Petersson metric. Under the
Weil-Petersson metric, $\mathcal M$ is a K\"ahler
orbifold~\footnote{A K\"ahler orbifold metric is a K\"ahler
metric on the smooth part of the orbifold and lifts to an
invariant K\"ahler metric on each local uniformization.
See~\cite{Ruan} for details of  orbifolds
and vector bundles over orbifolds.}.

Let $\Omega$ be a (nonzero) holomorphic $(n,0)$-form
on
$Z_t$. Define $\Omega\lrcorner\iota(\frac{\pa}{\pa
t_i})$ to be the contraction of $\Omega$ and
$\iota(\frac{\pa}{\pa t_i})$. The Weil-Petersson metric can
 be re-written as (cf. ~\cite{T1}):
\begin{equation}\label{tian}
\left(\frac{\pa}{\pa t_i},\frac{\pa}{\pa \bar t_j}
\right)_{WP}
=-\frac{\int_{Z_t}\Omega\lrcorner\iota(\frac{\pa}{\pa
t_i})\wedge\overline{\Omega\lrcorner\iota(\frac{\pa}{\pa
t_j})} }{\int_{Z_t}\Omega\wedge\bar\Omega}.
\end{equation}

The Weil-Petersson metric is the most natural metric on the moduli
space. Unfortunately, it does not have a very good curvature
property. In~\cite{Lu3}, another natural metric called the
Hodge metric was defined.
In this paper,  we  use the notations  in~\cite{LS-1}
for the Hodge metric.

Recall that for an $n$-dimensional compact complex
manifold $X$ with polarization, for any $0\leq k\leq n$, we have
the decomposition of the Hodge bundles
\[
H^k(X,\mathbb C)=\underset{p+q=k}{\oplus}
H^{p,q}(X,\mathbb C).
\]
By the Lefschetz decomposition theorem, we can further decompose
the Hodge bundles $H^{p,q}$ into its primitive parts as follows.

\newcommand{\C}{\mathbb C}

Define $L: H^k(X,\C)\rightarrow H^{k+2}(X,\C)$ by
$[\alpha]\rightarrow[\alpha\wedge\omega]$, where $\omega$ is the
curvature form of the ample line bundle over $X$. Define the
primitive cohomology group $P^k(X,\C)$ to be the kernel of
$L^{n-k+1}$ on $H^k(X,\C)$. Let $PH^{p,q}=P^k(X,\C)\cap
H^{p,q}(X)$.

The Lefschetz decomposition theorem states that
\begin{equation}\label{decomp}
H^{p,q}(X)=PH^{p,q}\oplus L(PH^{p-1,q-1})\oplus\cdots
\oplus L^r(PH^{p-r,q-r}),
\end{equation}
where $r={\rm Min}\,(p,q)$.

Define
\begin{equation}\label{qqq}
Q(\eta_1,\eta_2)=\int\eta_1\wedge\eta_2\wedge
\omega^{n-2k}.
\end{equation}
for $\eta_1,\eta_2\in H^k(X,\C)$. Then $Q$
extends to a bilinear form on $H^*(X,\C)$. The
Riemann-Hodge relations are
\begin{enumerate}
\item $Q(\eta_1,\eta_2)=0$, if $\eta_1\in
PH^{p_1,q_1}, \eta\in PH^{p_2,q_2}$, but
$p_1+p_2\neq q_1+q_2$;
\item $(\sqrt{-1})^{q-p}Q(\eta_1,\bar \eta_1)>0$
if $0\neq\eta_1\in PH^{p,q}$.
\end{enumerate}

The second Riemann-Hodge relation defines a Hermitian inner
product on the primitive harmonic $(p,q)$ forms:
\begin{equation}\label{qqq1}
<\eta_1,\eta_2>=(\sqrt{-1})^{q-p}Q(\eta_1,\bar \eta_1).
\end{equation}
When $X$ is a polarized Calabi-Yau manifold
$Z$, the above inner product is
equivalent to the inner product induced from the Ricci-flat
metric(cf. ~\cite{Wells}, ~\cite{FL-a}):

\begin{theorem}\label{Riem-Hodge}
%[The second Hodge-Riemann relation]
Let $\phi\in H^{p,q}(Z),p\geq
q$ and let
\[
\phi=\phi_0+L\phi_1+\cdots+L^q\phi_q
\]be the decomposition corresponding to (\ref{decomp}).
Then we have
\begin{align}\label{422}
\begin{split}
&
<\phi,\phi>\\
&=(-1)^{\frac 12(p+q)(p+q+1)}
\sum_{k=0}^{q}(-1)^k(n-p-q+2k)!\int_M||\phi_k||^2dV_{CY},
\end{split}
\end{align}where $||\cdot||$ is the metric induced from
the Ricci-flat
metric of $Z$.
\end{theorem}

\qed

If we make the relative version of the above settings, we get
bundles
$PR^q\pi_{*}\Omega_{\frak X/\mathcal M}^p\rightarrow \mathcal M$
in place of the cohomology groups. These bundles are called Hodge
bundles. They are in fact the vector bundles of K\"ahler
orbifolds (cf.~\cite{Ruan}).

The Kodaira-Spencer map $T_t\mathcal M\rightarrow
H^1(Z_t,\Theta_t)$ gives a bundle map
\[
\frac{\pa}{\pa t_i}: PR^q\pi_{*}\Omega^p_{\frak X/\mathcal M}\rightarrow
PR^k\pi_{*}(\C)/PR^q\pi_{*}\Omega^p_{\frak X/\mathcal M}
\]
for $k\leq n$
by differentiation. In this way, we have a natural bundle
 map (of vector bundles over \ka orbifold)
\begin{equation}\label{12-1}
T\mathcal M\rightarrow\underset{p+q=k}{\oplus}
{\rm Hom}\,(PR^q\pi_{*}\Omega^p_{\frak X/\mathcal M},
PR^k\pi_{*}(\C)/PR^q\pi_{*}\Omega^p_{\frak X/\mathcal M}).
\end{equation}

\begin{definition}\label{def22}
For each $t\in\mathcal M$ and $Z=Z_t$ with  the polarized
Ricci flat metric, Theorem~\ref{Riem-Hodge} defines  Hermitian
metrics on the bundles
$PR^q\pi_{*}\Omega^p_{\frak X/\mathcal M}\rightarrow\mathcal
M$. Let
$h_{P,k}$ be the  pull back of the natural Hermitian metric
on the bundle
$\underset{p+q=k}{\oplus}{\rm Hom}
(PR^q\pi_{*}\Omega^p_{\frak X/\mathcal M}\rightarrow
PR^k\pi_{*}(\C)/PR^{q}\pi_*\Omega^p_{\frak X/\mathcal M})$ to
$T\mathcal M$ for $k\leq n$. We use $\omega_{PH^k}$ to denote the
corresponding K\"ahler forms for $k\leq n$.\footnote{That is, if
$h_{H^k}=(h_{H^k})_{i\bar j} dt_i\otimes d\bar t_j$, then
$\omega_{H^k}=\frac{\sqrt{-1}}{2\pi}(h_{H^k})_{i\bar j}
dt_i\wedge d\bar t_j$.}   According to (\ref{decomp}), we
define
\begin{equation}\label{12-2}
\omega_{H^k}=\omega_{PH^k}+\omega_{PH^{k-2}}+\cdots.
\end{equation}
We call both $\omega_{H^k}$ and $\omega_{PH^k}$ to be the
generalized Hodge metrics.
\end{definition}

\begin{remark}\label{rk1}
The generalized Hodge metric is a generalization of the
Hodge metric defined by the second author~\cite{Lu5}. In fact,
 it is proved in~\cite{LS-1} that
$$\omega_{PH^n}=\omega_H,$$ the latter being the Hodge metric.
\end{remark}

Because of the possible degeneration of the action~(\ref{12-1}), the
generalized Hodge metric is only positive semi-definite; hence,
it was only
a pseudo-metric. However, it enjoys similar ``curvature"
properties
of the Hodge metric proved by the second author~\cite{Lu3}.
The generalized Hodge metrics are bounded by the  Poincar\'e
metric. See Appendix A for more details.

We  have the following relations between the Hodge metric and
the Weil-Petersson metric:

\begin{theorem}\label{thm22}
We use the above notations. The Hodge metric
$\omega_H$ defines a K\"ahler metric (i.e. $d\omega_H=0$).
Furthermore, The bisectional curvature of $h$ is non-positive and
the holomorphic sectional curvature
and the Ricci curvature are negative away from zero. In
particular, we have
\begin{enumerate}
\item If $n=2$, then $\omega_H=2\omega_{WP}$;
\item If $n=3$, then $\omega_H=(m+3)\omega_{WP}
+{\rm Ric}(\omega_{WP})$;~\cite{Lu5}
\item If $n=4$, then $\omega_H=(2m+4)\omega_{WP}
+2{\rm Ric}(\omega_{WP})$,~\cite{LS-1}
\end{enumerate}
where $n=\dim Z$ and $m=\dim \mathcal M$.
\end{theorem}

\begin{definition}
We call a Calabi-Yau manifold primitive if $\omega _{H^{k}}=0$ for all $k<n$.
\end{definition}

A direct consequence of the above definition,
equation (\ref{12-2})  and Remark~\ref{rk1}
is that $\omega_{PH^n}=\omega_{H^n}=\omega_{H}$
when $Z$ is primitive. Thus, the only
interesting generalized Hodge metric is the original Hodge metric.
The class of primitive
Calabi-Yau manifolds includes many important examples.
\begin{example}
A Calabi-Yau three-fold is by definition simply
connected.
It is easy to see that the actions of $T\mathcal{M}$ on lower
degree Hodge bundles are trivial; hence, it is primitive.
Similarly, a Calabi-Yau four-fold with vanishing $h^{1,2}$ is also
primitive.
\end{example}

\begin{example}
More generally, due to the hard Lefschetz theorem, all the Calabi-Yau
hyper-surfaces and complete intersections in projective spaces are primitive.
\end{example}

We give the explicit formulae for the generalized
Hodge metrics in the next proposition:
\begin{proposition}
\label{prop1}
Let $c_1(E)$ be the Ricci form  of a vector bundle $E$. Then we have
\begin{equation}\label{five}
\omega _{PH^{k}}=\sum_{0\leq p\leq
k}pc_{1}(PR^{k-p}\pi_{*}\Omega^{p}_{\frak X/ \mathcal M}).
\end{equation}
\begin{equation}\label{six}
\omega _{H^{k}}=\sum_{0\leq p\leq
k}pc_{1}(R^{k-p}\pi_{*}\Omega^{p}_{\frak X/ \mathcal M}).
\end{equation}
for $k\leq n$.

\end{proposition}

{\bf Proof.} Fixing a $k\leq n$, we define the Hodge bundles
$\mathcal F^k_k,\cdots,\mathcal F^0_k$ to be
\[
\mathcal F_k^p=PR^{0}\pi_*\Omega^{k}_{\frak X/\mathcal M}\oplus
\cdots\oplus PR^{k-p}\pi_*\Omega^{p}_{\frak X/\mathcal M}
\]
for $p=0,\cdots,k$. Thus  for $q=k-p$,
\[
PR^{q}\pi_*\Omega^{p}_{\frak X/\mathcal M}=\mathcal F_k^p/
\mathcal F_k^{p+1}.
\]
In terms of the curvatures, we have

\begin{equation}\label{q0}
 c_1(PR^{q}\pi_*\Omega^{p}_{\frak X/\mathcal M})=c_1(\mathcal
F_k^p) -c_1(F_k^{p+1}).
\end{equation}

By the Abel summation formula, we have
\begin{equation}\label{q1}
\sum_{0\leq p\leq k}pc_{1}(PR^{k-p}\pi_{*}\Omega^{p}_{\mathfrak X/
\mathcal M})=c_1(\mathcal F_k^k)+\cdots+c_1(\mathcal F_k^1)
+c_1(\mathcal F_k^0).
\end{equation}

Each $\mathcal F_k^p$ is a sub-bundle of the flat bundle $\mathcal
F_k^0 =PR^k\pi_*\C$. Let $t_1,\cdots,t_m$ be the local holomorphic
coordinate of $\mathcal M$ and let the bundle map
\[
\frac{\pa}{\pa t_i}: \mathcal F_k^p\rightarrow
 \mathcal F_k^0/\mathcal F_k^p, \,1\leq i\leq m
\]
be represented by the matrix
\[
\frac{\pa\Omega_\alpha}{\pa t_k}=b_{k\alpha\mu}T_\mu,
\]
where $\Omega_\alpha$ and $T_\mu$ are the basis of $\mathcal
F_k^p$ and $\mathcal F_k^0/\mathcal F_k^p$, respectively. Then the
first Chern class can be represented by
\begin{equation}\label{12-3}
c_1(\mathcal F_k^p)=\frac{\sqrt{-1}}{2\pi}
\sum_{\alpha,\mu}b_{k\alpha\mu}\bar b_{l\alpha\mu}dt_k\wedge d\bar
t_l
\end{equation}
for $0\leq p\leq k$. (\ref{five}) follows from the definition of
$\omega_{PH^k}$. (\ref{six}) follows from (\ref{five}) and
(\ref{decomp}). The proof is completed.

\qed

The curvature computation is a natural generalization of the
similar result in \cite{Gr},  where only the middle
dimensional primitive Hodge structure was considered.

\begin{remark}\label{12-7}
The Weil-Petersson metric is the curvature
of the first Hodge bundle:
$$ \omega_{WP}=c_1(R^0\pi_*(\Omega^n_{\mathfrak X/\mathcal
M}))$$
by~\eqref{tian},~\eqref{12-3}.
Thus by the above equation, Remark~\ref{rk1}, and
Proposition~\ref{prop1},  the Weil-Petersson metric, the Hodge
metric and the  generalized Hodge metrics are
the Ricci curvatures of the combination of the Hodge bundles.
\end{remark}

With the above
interpretation of the metrics, the following is
obvious:

\begin{corollary}
\label{coro1}
Using the above notations, for $n\geq 2$, we have
\begin{equation}
\omega _{H}\geq 2\omega _{\text{WP}}.
\end{equation}
\end{corollary}

{\bf Proof.} By Remark~\ref{12-7} and Serre Duality, we get
\begin{equation}
\omega _{_{WP}}=c_{1}(R^0\pi_*(\Omega^n_{\mathfrak X/\mathcal
M}))=-c_1(R^n\pi_*(\mathcal O)). \label{1}
\end{equation}

By (\ref{q0}) and the fact that $\mathcal F^0_n$ is flat, we have
\begin{equation}\label{12-5}
c_1(\mathcal F_n^1)=-c_1(R^n\pi_*(\mathcal O))=\omega _{_{WP}}.
\end{equation}

Also, since $\mathcal F_n^{n+1}=0$,
\begin{equation}\label{12-6}
c_1(\mathcal F^n_n)=c_1(R^0\pi_*(\Omega^n_{\mathfrak X/\mathcal
M}))=\omega_{WP}.
\end{equation}

According to (\ref{12-3}), $c_1(\mathcal F_n^p)\geq 0$ for all $p$. Hence, when $n\geq 2$,
by  (\ref{five}), (\ref{q1}), (\ref{12-5}) and
(\ref{12-6}), we have
\begin{equation}
\omega _{_H}\geq c_1(\mathcal F_n^n)+c_1(\mathcal F_n^{1})=2\omega _{WP}.
\end{equation}
The proof is finished.

\qed

\section{BCOV Torsion}\label{section3}

\setcounter{equation}{0} \setcounter{theorem}{0}BCOV torsion was
first defined by Bershadsky-Ceccotti-Ooguri-Vafa in their study of
the Physics Mirror Symmetry.\ It was constructed as the
partition function for the N=2 SCFT. In their breakthrough works
\cite{BCOV2,BCOV1}, the torsion was determined using the Mirror
Symmetry Conjecture.\ One amazing consequence is that, given the
local expansion of the BCOV torsion, they were able to give a
prediction of counting numbers of embedded elliptic curves  of all degrees
in a given Calabi-Yau manifold. The prediction matches
all known low degree cases. Furthermore, they also discussed the
higher genus cases based on genus one computation. Notice that a
similar prediction of counting rational curves for quintics,
which invoked
intensive mathematical research, was
first made by Candelas et. al.\cite{COGP}.  Through the
fundamental works of Kontsevich,\ Givental,\ Lian-Liu-Yau and
many others (see~\cite{CoxKatz} for a complete reference), the
prediction has been mathematically verified. It is thus of
crucial interest to understand the BCOV torsion in terms of
algebraic and differential geometry.

First, we give the following

\begin{definition}
The BCOV torsion of a Calabi-Yau manifold is
\begin{equation}
T=\prod_{1\leq p,q\leq n}(\det \Delta'
_{p,q})^{(-1)^{p+q}pq},
\end{equation}
where $\Delta _{p,q}$ is the $\overline{\partial }$-Laplace operator
on $(p,q)$ forms
with
respect to the Ricci-flat metric on a fiber;
$\Delta _{p,q}^{\prime }$ represents the non-singular part of
$\Delta_{p,q}$; the determinant
is taken in the sense of zeta function regularization.
\end{definition}

The BCOV torsion is an analytic torsion  in
the sense of Ray
and Singer~\cite{RS-1}. To see this, we define a holomorphic
coefficient vector bundle over a polarized Calabi-Yau
manifold $Z$,
\begin{equation}
E=\bigoplus\limits_{p=1}^{n}(-1)^{p}p\Omega
^{p}(\frak{X}/\mathcal{M)}. \label{3-1}
\end{equation}
$E$ inherits a natural Hermitian metric induced from the metric on the
relative tangent bundle. According to \cite{KnM}, there exists a
corresponding determinant line bundle over $\mathcal{M}$, which is
defined to be

\begin{equation}
\lambda =\bigwedge_{0\leq p,q\leq n}(\det (H^{p,q}(Z,E,\overline{\partial }%
))^{(-1)^{p+q}p},  \label{3-2}
\end{equation}
where $H^{p,q}(Z,E,\overline{\partial })=R^{q}\pi _{\ast }\Omega ^{p}(\frak{X%
}/\mathcal{M)}$  are holomorphic vector bundles over
$\mathcal{M}$; and we identify the cohomology groups with the
corresponding harmonic forms with respect to the natural induced
metrics on various spaces.

There are two natural metrics defined on $\lambda$. The usual L$^{2}$
metric is defined by the harmonic
forms, and the Quillen metric is given by

\begin{equation}\label{3-3}
\left\| \cdot \right\| _{Q}=\left\| \cdot \right\| _{L^{2}}T.
\end{equation}
Following the comprehensive studies on the determinant line
bundles and the associated Quillen metrics in
\cite{BGS1,BGS2,BGS3}, where local index techniques
are extensively used, people have made many progress in this
direction with  applications in many branches of
mathematics. See  \cite{B-ICM}  for more references and survey.

The formulation that we have in (\ref{3-1}), (\ref{3-2}) and
(\ref{3-3}) make it possible to study BCOV torsion in the
framework developed in~\cite{BGS1,BGS2,BGS3}. First, we verify the
following:

\begin{lemma}\label{localK}
$\frak X$ satisfies the local K\"ahlerian property defined
in~\cite{BGS1} and~\cite{BGS2}.
\end{lemma}

\textbf{Proof.} Assume that $U$ is a coordinate open set of $\mathcal{M}$.%
\footnote{%
If $U$ contains a quotient singularity of $\mathcal{M}$. Then $\mathfrak{X}%
|_{U}$ is K\"{a}hlerian in the following sense: let
$(\hat{U},p,U)$ be the local uniformization system. That is,
$\hat{U}$ is a finite cover of $U$ and $\hat{U}$ is homeomorphic
to the unit ball of $\mathbb{C}^{m}$. $U=\hat{U}/G$ by a finite
group $G$. Let $\hat{\mathfrak{X}}_{\hat{U}}$ be the pull back of
$\mathfrak{X}|_{U}$ through the map $p:\hat{U}\rightarrow U$. Then
we say $\mathfrak{X}|_{U}$ is K\"{a}hlerian if
$\mathfrak{\hat{X}}|_{\hat{U}}$ is and the K\"{a}hler metric is
$G$-invariant. The proof of the lemma will go through in this
case.} Let $h_t$ be the Hermitian metric of $L$ over the fiber
$\pi^{-1}(t)$, where $L$ is the polarization of $Z$. The form
$\frac{\sqrt{-1}}{2\pi} \partial\bar\partial\log h_t$, when
restricted on $\pi^{-1}(t)$, is positive. Let $t_1,\cdots, t_m$ be
local holomorphic coordinate system of $U$. Then
\[
\frac{\sqrt{-1}}{2\pi}a\sum_i dt_i\wedge d\bar t_i
+\frac{\sqrt{-1}}{2\pi} \partial\bar\partial\log h_t
\]
gives the local K\"ahler metrics on $\pi^{-1}(U)$ for
$a\in\mathbb R$ large enough.

\qed

The following Lemma of
Beshadsky-Cecotti-Ooguri-Vafa~\cite{BCOV2,BCOV1} characterizes
the geometry of the coefficient bundle $E$.

\begin{lemma}\label{34}
For $E$ associated with the induced metric from that of
$\frak{X},$ as forms on $\frak{X}$ we have
\begin{equation}
Td(T^{1,0}\frak{X/}\mathcal{M})ch(E)=-c_{n-1}+\frac{n}{2}c_{n}-\frac{1}{12}%
c_{1}c_{n},
\end{equation}
where $c_{i}=c_{i}(T^{1,0}\frak{X/}\mathcal{M})$.
\end{lemma}

\qed

As a consequence,

\begin{proposition}\label{prop33}
The Quillen metric is the potential of
$\chi_Z/12$ times the Weil-Petersson metric.
\begin{equation}
c_{1}(\lambda ,\left\| \cdot \right\| _{Q})=\frac{\chi _{Z}}{12}\omega _{%
\text{WP}}.
\end{equation}
\end{proposition}

{\bf Proof.}
This also appears in~\cite{BCOV1}. We include the
proof  here for
completeness. First  notice that, by using the Ricci-flat metric on the
relative tangent bundle and  Remark \ref{rk1}, we have
\begin{equation}\label{3A}
c_{1}(T^{1,0}\frak{X/}\mathcal{M})=-c_1(H^{n,0})=-\omega _{WP}.
\end{equation}
Then, by Lemma~\ref{localK}, a direct application of the family
Grodenthick-Riemann-Roch theorem proved by
Bismut-Gillet-Soul\'e\cite{BGS1}, shows that, for $E$
in~\eqref{3-1},
\begin{equation}
c_{1}(\lambda ,\left\| \cdot \right\| _{Q})=
[\int_{Z}\text{Td}(T^{1,0}\mathfrak X/\mathcal M)\text{ch}%
(E)]^{(1,1)}.
\end{equation}
By Lemma~\ref{34},~\eqref{3A} and the Gauss-Bonnet formula,
we have :
\begin{equation}
c_{1}(\lambda ,\left\| \cdot \right\| _{Q})=
[\int_{Z}{1\over 12}
\omega_{WP}c_n(T^{1,0}\frak{X/}\mathcal{M})]^{(1,1)}=\frac{\chi
_{Z}}{12}\omega _{%
\text{WP}}.
\end{equation}
\qed

The next proposition relates the curvature of $\lambda $ with
respect to the L$^{2}$-metrics to the generalized Hodge metrics.

\begin{proposition}\label{prop35}
Using the notations in the above and in the previous section,
we have
\begin{equation}
c_{1}(\lambda ,\left\| \cdot \right\|
_{L^{2}})=\sum\limits_{i=1}^{n}(-1)^{i}\omega _{H^{i}}.
\end{equation}
\end{proposition}

{\bf Proof.}
This is due to Proposition \ref{prop1}. It is an easy computation
to show that
\begin{equation}
c_{1}(\lambda ,\left\| \cdot \right\|
_{L^{2}})=
\sum_{0\leq p,q\leq n}(-1)^{p+q}pc_{1}(R^{q}\pi _{\ast }\Omega
^{p}(\frak{X}/\mathcal{M)})=\sum_{k=1}^{n}(-1)^{k}\omega _{H^{k}}.
\end{equation}
The proof is complete.

\qed

{\bf Proof of Theorem~\ref{main}.} Using the relation of the
$L^2$ metric, the Quillen metric and the BCOV torsion,
~\eqref{markertest} follows from
~\eqref{3-3}, Proposition~\ref{prop33},
and Proposition~\ref{prop35}.
The equation~\eqref{marker2} follows from
~\eqref{markertest} and Remark~\ref{rk1}.

\qed

Theorem~\ref{main} is the  explicit relation between two
kinds
canonically defined K\"{a}hler metrics: the
Weil-Petersson metric and the generalized Hodge metrics,
on the moduli space. The
surprising fact is that the bridge is the BCOV torsion, a
spectral invariant of Ricci-flat metrics, which is also of
its own significance in physics literature.

{\bf Proof of Corollary~\ref{cor14}.} The global
Poincar\'e
metric defined at the beginning of Appendix A is asymptotic to
the  Poincar\'e metric in ~\eqref{eq15}. So they are
equivalent on $(\Delta^*)^l\times\Delta^{m-l}$.  The corollary
is thus follows from Theorem~\ref{main},
Corollary~\ref{coro1} and Theorem~\ref{A1}.

\qed

Asymptotic expansion of $T$ near points of maximal
monodromy degeneration will give the prediction for counting
elliptic curves in the Calabi-Yau manifolds. The above result
determines the asymptotic behavior of the BCOV torsion up to a
(possibly multi-valued) pluri-harmonic function.

\smallskip

{\bf Proof of Corollary~\ref{mo}.}
If $N$ is a smooth submanifold in the smooth part of
 $\mathcal
M$, the corollary follows from the ordinary Stokes theorem. In
general, we  generalize the Stokes theorem into the singular
case. Let's first define  the Hodge and the
Weil-Petersson volumes  on $N$.

Let $\mu$ be a smooth $2k$ form of the orbifold $\mathcal M$,
where
$k=\dim \,N$. Then
we can  define a measure on $N$ as follows:

Let $x\in N\subset \mathcal M$ and let $(\hat U, p, U)$
be a local uniformization of $\mathcal M$ at $x$; i.e.,
$p: \hat U\to U=\hat U/G$ for a finite group $G$. Let $V=U\cap N$.
Let $\hat N\subset \hat U$ be the pre-image of $N$ under $p$.
For $y\in V$, define $n(y)$ to be the multiplicity
of the map $\hat N\rightarrow N$. We can define a measure
$\mu$
on $V$ by
\[
\mu(V)=\int_{\hat V}{\frac {1}{n(y)}} p^*(\mu)|_{\hat N},
\]
where $\hat V$ is the pre-image of $V$ under $p$. The
generalized
Hodge metrics, the Weil-Petersson metric, and
$\frac{\sqrt{-1}}{2\pi}\pa\bar\pa\log \,T$,
through this definition,
 define the corresponding measures on $N$.
In particular, the Hodge and the Weil-Petersson volumes are
defined.

The function $n(y)$ defined above is in fact independent of the
choice of the local uniformization so it is a global function
on
$N$. It is a
constant on some Zariski open set
$N'$
of $N$. Because of this, it is sufficient to prove that
\begin{equation}
\text{Vol}_{H}(N')=[\frac{(-1)^{n}}{12}\chi
_{Z}]^{k}\text{Vol}_{{\normalsize  WP}}(N').
\end{equation}
Thus the corollary follows from the following version of
 Stokes theorem:

\begin{lemma}\label{lem36}
Let $\eta$ be a smooth $2k-1$ form on the
orbifold $\mathcal M$. Then
$d\eta$ defines a measure $\mu(\eta)$ on $N$. We have
\[
\int_{N_{{\rm reg}}}\mu(\eta)=0,
\]
where $N_{reg}$ is the smooth part of $N$.
\end{lemma}

{\bf Proof.}
Let $N''=N-N_{reg}\cap N'$. Let $\rho$ be a smooth function
such that:
(1). $\rho=0$ if $dist(x,N'')<\eps$; (2). $\rho=1$
if $dist(x,N'')>2\eps$; and $|\nabla\rho|\leq
3/\eps$.

By the ordinary Stokes Theorem,
\[
\int_Nd(\rho\eta)=0.
\]
Thus we have
\[
\int_{N_{reg}\cap N'}\rho d\eta+\int_Nd\rho\wedge\eta=0.
\]
Since $N''$ is compact, its Hausdorff measure is finite. Thus
\[
\left|\int_Nd\rho\wedge\eta\right|
\leq C \cdot\frac{3}{\eps} \cdot\pi\eps^2
(\text{H}(N''))\rightarrow 0.
\]
Let $n(y)=const$ on $N'$. Then
\[
0\leftarrow const\cdot\int_{N_{\rm reg}\cap N'}\rho d\eta
\rightarrow \int_N\mu(\eta)
\]
This completes the proof.

\qed

{\bf Proof of Corollary~\ref{mi}.}
Suppose $N$ is a complete curve of $\mathcal M$. Then
by Corollary~\ref{mo},
\[
\text{Vol}_H(N)=\frac{(-1)^n}{12}\chi_Z\text{Vol}_{WP}(N).
\]
By the assumption, we thus have
\[
\text{Vol}_H(N)<2\text{Vol}_{WP}(N),
\]
unless $N$ is  of dimension $0$. However, the
above inequality contradicts  to Corollary~\ref{coro1}.

\qed

\begin{remark}
All the Calabi-Yau three-folds are primitive.\ By the Physics
Mirror Symmetry, it is conjectured that all the Calabi-Yau
three-folds occur in pairs (mirror pair); and the Euler characteristics of paired
Calabi-Yau three-folds are differed only by signs.\ Hence,
Corollary~\ref{mi} claims that ``more than half'' of moduli of
polarized Calabi-Yau three-folds contains no complete
subvariety.
\end{remark}

\newcommand{\bb}{\frac{\ii}{2\pi}}

We proceed to discuss the obstruction to the existence of complete curves in
the Calabi-Yau moduli, in some more specific situations.

\begin{remark}
For  a primitive Calabi-Yau four-fold $Z$,  let
$N$ be a complete curve in the moduli space $\mathcal M$.
By Theorem~\ref{thm22}, we have
\begin{equation}\label{b1}
\int_N\omega_H=(2m+4)\int_N\omega_{WP}+2\int_N{\rm Ric}
(\omega_{WP}).
\end{equation}
By Theorem~\ref{main}, we have
\begin{equation}\label{b2}
\int_N\omega_H=\frac{\chi_Z}{12}\int_N\omega_{WP}.
\end{equation}
If $\chi_Z>24(m+2)$, from ~\eqref{b1}
and~\eqref{b2}, we have
\begin{equation}\label{b3}
\int_N{\rm Ric}(\omega_{WP})\geq 0.
\end{equation}
On the other side, the Ricci curvature of $\omega_H$
is negative, so
\[
\int_{N}{\rm Ric}(\omega_H)<0.
\]
The above inequality contradicts to ~\eqref{b3} because
${\rm Ric}(\omega_H)-{\rm
Ric}(\omega_{WP})=-\frac{\sqrt{-1}}{2\pi}
\pa\bar\pa\log(\frac{\omega_H^m}{\omega_{WP}^m})$, and thus
using Lemma~\ref{lem36}, the integration of the two
Ricci curvatures are the same. This is a contradiction.
Thus if $\chi_Z>24(m+2)$, there is no complete curve in
$\mathcal M$. The similar result in the case
of Calabi-Yau threefold is trivial and the similar results
for high dimensional Calabi-Yau manifolds are still unknown.

\end{remark}

\begin{remark} When  the dimension of $\mathcal M$ is
2, if there
exists a complete curve $C\subset{\mathcal M}$, then
\begin{equation}\label{poi}
\int_C{Ric_{WP}}=\chi_C+C.C,
\end{equation}
where
$C.C$ is the self-intersection number. Let
$\omega'=\omega_{WP}|_C$.
In local coordinate, let
\[
\omega'=\bb\,h_{1\bar 1} dt_1\wedge d\bar t_1,
\]
and let
\[
\omega_{WP}=\bb\,\sum_{i,j=1}^2h_{i\bar j}
dt_i\wedge d\bar
t_j.
\]
Assume that at a point $x$ $h_{i\bar j}=\delta_{ij}$. Then we
have
\[
\pa_1\bar\pa_1\log h_{1\bar 1}=\frac{\pa^2 h_{1\bar 1}}{\pa
t_1\pa\bar t_1}-h_{1\bar 1}^{-1}\left|\frac{\pa h_{1\bar
1}}{\pa t_1}\right|^2.
\]
Since
\[
h_{1\bar 1}^{-1}\left|\frac{\pa h_{1\bar
1}}{\pa t_1}\right|^2\leq
h^{i\bar j}\frac{\pa h_{1\bar
j}}{\pa t_1}\cdot\frac{\pa h_{i\bar 1}}{\pa\bar t_1},
\]
we have
\[
\pa_1\bar\pa_1\log h_{1\bar 1}\geq
R_{1\bar 11\bar 1},
\]
where $R_{i\bar jk\bar l}$ is the curvature tensor of
$\omega_{WP}$. By the (generalized) Strominger formula
(cf.~\cite[Theorem 3.1]{LS-1}, ~\cite{S}, and ~\cite{Wang1}),
we have
\[
R_{1\bar 12\bar 2}\leq 1.
\]
Thus we have the following
\[
{\rm Ric}(\omega')
\leq {\rm Ric}(\omega_{WP})+1.
\]
Comparing the above equation with~\eqref{poi}, we have

\[
C.C\geq -1.
\]
 This would be a new
obstruction to the existence of complete curves in these cases.

\end{remark}
%\begin{corollary}
%Assume the Calabi-Yau manifold $Z$ is primitive.\ $%
%(-1)^{n+1}\chi _{Z}>-24$. Let $\mathcal{\bar M=M\cup\partial M}$
%be the compactified muduli, with a given hyperplane section $L$.
%Then\\
%either $\mathcal{\partial M.\partial M}.L^{n-2}<0$, or\\
%the Poncar\'e dual of $\mathcal {\partial M}$ generates $H^2(M,C)$.
%\end{corollary}

Notice that the results we proved  above are only on the
smooth moduli. It is very interesting to see the
corresponding results for the compactified moduli. This kind
of global results will be achieved by finer local analysis of
the BCOV torsion near the boundary of the moduli.

\appendix
\section{Poincar\'e Metric and Generalized Hodge Metrics}

\setcounter{equation}{0} \setcounter{theorem}{0}
In this Appendix, we prove that the generalized Hodge metrics
are bounded by the Poincar\'e metric. The main
 result of this Appendix, Theorem~\ref{A1}, is a degenerate
version of Yau's Schwarz Lemma (cf.~\cite{Y3}) from \ka
manifolds to Hermitian manifolds.

Calabi-Yau moduli $\mathcal M$ is
quasi-projective. The smooth part $\mathcal M_{\rm reg}$
of $\mathcal M$ allows a compactification
$\bar{\mathcal M}$, where $\bar{\mathcal M}$ is a compact smooth manifold
such
that $\bar{\mathcal M}\backslash{\mathcal M}_{\rm reg}$ is a
divisor of normal
crossing. For the pair of manifolds
$(\bar{\mathcal M},\mathcal M_{\rm reg})$, we define
a \ka metric $\omega_{GP}$ called the global Poincar\'e
metric as follows: in a neighborhood
$(\Delta^*)^l\times\Delta^{m-l}$ of
$x\in\bar{\mathcal M}\backslash{\mathcal M}$,
 the global Poincar\'e metric $\omega_{GP}$
is  asymptotically the Poincar\'e metric $\omega_P$ defined
in~\eqref{eq15}:
\[
\omega_{GP}\sim \omega_P=\sum_{i=1}^l\ii\frac{1}{
|z_i|^2(\log\frac{1}{|z_i|})^2} dz_i\wedge d\bar z_i
+\sum_{i=l+1}^m
\ii dz_i\wedge d\bar z_i.
\]
The global Poincar\'e metric $\omega_{GP}$ is a complete \ka
metric on $\mathcal M$ whose Ricci curvature is  bounded from
below. For the detailed construction of $\omega_{GP}$, see
~\cite{LS-2}.

The main result of this section is the following:
\begin{theorem}\label{A1}
Using the above notations, we have
\[
\omega_{PH^k}\leq C\omega_{GP}
\]
for $0\leq k\leq n$ and  a constant $C$ depending on  the
lower bound of
\,${\rm Ric}(\omega_{GP})$, dimension of the Calabi-Yau
manifolds and dimension of the moduli space $\mathcal M$.
\end{theorem}

We use $H^{p,q}$ to denote the bundle $PR^q\pi_*\Omega^p_{
\mathfrak
X/\mathcal M}$ for $0\leq p,q\leq n$.  Let
$F_k^p=H^{k,0}\oplus\cdots\oplus H^{p,k-p}$ for
$p=0,\cdots,k$. We also assume that
$H^{k+2,-2}=H^{k+1,-1}=H^{-1,k+1}=H^{-2,k+2}=0$ and $F_k^{k+1}=0$,
$F_k^{-1}=F^0$, for the sake of simplicity.

First we try to express the generalized Hodge metrics in local
coordinates. Fix $k\leq n$, $p\leq k$ and $q=k-p$. Let
$\{\Omega_{p,i}\}$, $i=1,\cdots, h^{p,q}$ be a local
frame of $H^{p,q}$.

\begin{definition}
Let $(t_1,\cdots,t_m)$ be a holomorphic local coordinate at
a point of $\mathcal  M_{\rm reg}$. We define
$D_\alpha\Omega_{p,i}\in H^{p-1,q+1}$
to be the projection of $\pa_\alpha\Omega_{p,i}
=\frac{\pa}{\pa t_\alpha}\Omega_{p,i}$
to $H^{p-1,q+1}$ with respect to the bilinear form
$Q(\,\,,\,\,)$ in~\eqref{qqq}.
\end{definition}

For simplicity, we shall use $(\,\,,\,\,)$ in stead of the
bilinear form $Q$ in~\eqref{qqq}.
With the above notation,
\begin{equation}\label{first}
(g_p)_{i\bar j}=<\Omega_{p,i}, \overline{\Omega_{p,j}}>
=(\ii)^{q-p}(\Omega_{p,i}, \overline{\Omega_{p,j}})
\end{equation}
 is the Hermitian metric matrix of $H^{p,q}$
for $p=0,\cdots,k$. It is thus easy
to see that

\begin{proposition}
The generalized Hodge metric matrix for the local
coordinate system $(t_1,\cdots,t_m)$ with respect to
$PH^{k}$, defined in Definition~\ref{def22}, is
\begin{equation}\label{second}
h_{\alpha\bar \beta}=\sum_{p=0}^k(\ii)^{q-p+2}
g_p^{i\bar
j} (D_\alpha\Omega_{p,i},\overline{D_\beta\Omega_{p,j}}),
\end{equation}
where $(g_p^{i\bar j})$ is the inverse of $(g_p)_{i\bar j}$.
\end{proposition}

\qed

We proceed with two technical lemmas.
\begin{lemma}\label{lemmaA1}
\begin{equation}\label{third}
\bar\pa_\beta D_\alpha\Omega_{p,i}= g_p^{l\bar j}
<\bar\pa_\beta D_\alpha\Omega_{p,i},
\overline{\Omega_{p,j}}>\Omega_{p,l}.
\end{equation}
\end{lemma}

{\bf Proof.} We first claim that $\bar\pa D_\alpha \Omega_{p,i}\in
H^{p,q}$. To see this, let $\Omega_1\in F_k^{p+1}$. Then
\begin{equation}\label{prq}
(\bar\pa D_\alpha\Omega_{p,i}, \bar\Omega_1)
=-(D_\alpha\Omega_{p,i},\overline{\pa\Omega_1})=0.
\end{equation}
On the other hand, if $\Omega_2\in H^{p-1,q+1}$, we have the
decomposition
\[
\pa_\alpha\Omega_{p,i}=D_\alpha\Omega_{p,i}+B
\]
for $B\in H^{p,q}$. Furthermore,
\[
\bar\pa B\in F^p_k.
\]
Thus,
\begin{equation}\label{prq1}
\bar\pa D_\alpha\Omega_{p,i}\in F_k^p.
\end{equation}
Combining~\eqref{prq},~\eqref{prq1}, we have
$\bar\pa D_\alpha\Omega_{p,i}\in H^{p,q}$.
Writing $\bar\pa D_\alpha\Omega_{p,i}$ as the
linear combination of $\Omega_{p,i}$'s, we get \eqref{third}.

\qed

\begin{lemma}\label{lemmaA2}
If $A$ is a local section of $H^{p,q}$ and $B$ is a local
section
of $H^{p-1,q+1}$, then
\[
<D_\alpha A,\bar B>=<A,\overline{\bar\pa_\alpha B}>.
\]
\end{lemma}

{\bf Proof.} This follows from a straightforward computation:
\begin{align*}
&
<D_\alpha A,\bar B>=(\ii)^{q-p+2}(\pa_\alpha A,\bar B)\\
&
=-(\ii)^{q-p+2}(A,\overline{\bar\pa_\alpha B})
=<A,\overline{\bar\pa_\alpha B}>.
\end{align*}

\qed

{\bf Proof of Theorem~\ref{A1}.}
The generalized Hodge metrics are only semi-%
positive definite but not positive definite. If they were
positive definite, then using the similar method as
in~\cite{Lu3}, we should have been able to prove that the
holomorphic sectional curvatures of the metrics were
negative and bounded away from zero, and the holomorphic
bisectional curvatures of the metrics  were nonpositive. Thus
we could have used Yau's Schwarz Lemma~\cite{Y3} to get the
conclusion. The contribution of this Appendix is that we prove
the same result even if the generalized Hodge metrics fail to
be positive definite.

We assume all the notations in the previous sections. Write
the global Poncar\'e metric in local coordinates as:
\[
\omega_{GP}=\frac{\ii}{2\pi}\tau_{\alpha\bar\beta} d
t_\alpha\wedge d\bar t_\beta.
\]
Let $-C_1$ be the lower
bound of the Ricci curvature of $\omega_{GP}$ for
some constant $C_1>0$.
We define a smooth function $$f=\sum \tau^{\alpha\bar\beta}
h_{\alpha\bar\beta}$$
on $\mathcal M$. $f$ is nonnegative. If $f$ is bounded, then
$\omega_{PH^k}$ is bounded by $\omega_{GP}$.

For the rest of the Appendix,
we assume that at the given point $x$, $\tau_{i\bar
j}=\delta_{ij}$ and $d\tau_{i\bar j}=0$. Then at $x$,

\begin{equation}\label{impot}
\Delta f\geq -C_1 f+{\frac{\pa^2} {\pa t_\gamma\pa\bar
t_\gamma}}{h_{\alpha\bar\alpha}},
\end{equation}
where $\Delta$ is the Laplacian on $\mathcal M$ with respect to
$\omega_{GP}$.

We assume that at the point $x$, the frames $\Omega_{p,i}$
are chosen so that
$(g_p)_{i\bar j}=\delta _{ij}$, and
$\frac{\pa }{\pa t_\alpha}(g_p)_{i\bar j}=0$ for $p=0,\cdots,
k$ and $\alpha=1,\cdots, m$. Then a straightforward
computation gives
\begin{align}\label{fourth}
\begin{split}
&\frac{\pa^2 }{\pa t_\gamma\pa \bar
t_\gamma}h_{\alpha\bar\alpha}
=\sum_p(\ii)^{q-p+2}(-(R_p)_{j\bar i\gamma\bar\gamma})
(D_\alpha\Omega_{p,i},\overline{D_\alpha\Omega_{p,j}})\\
&+\sum_p(\ii)^{q-p+2}(\pa_\gamma D_\alpha\Omega_{p,i},
\overline{\pa_\gamma D_\alpha\Omega_{p,i}})\\
&+\sum_p(\ii)^{q-p+2}(\bar\pa_\gamma
D_\alpha\Omega_{p,i},
\overline{\bar \pa_\gamma D_\alpha\Omega_{p,i}})\\
&+\sum_p(\ii)^{q-p+2}(\bar\pa_\gamma\pa_\gamma
D_\alpha\Omega_{p,i},
\overline{D_\alpha\Omega_{p,i}})\\
&+\sum_p(\ii)^{q-p+2}(D_\alpha\Omega_{p,i},
\overline{\bar\pa_\gamma\pa_\gamma
D_\alpha\Omega_{p,i}}),
\end{split}
\end{align}
where $(R_p)_{j\bar i\gamma\bar\gamma}$ is the curvature
tensor of $g_p$ for $p=0,\cdots,k$.
By Lemma~\ref{lemmaA1}, we have
\begin{align}
&(\bar\pa_\gamma\pa_\gamma
D_\alpha\Omega_{p,i},
\overline{D_\alpha\Omega_{p,i}})=-(\bar\pa_\gamma
D_\alpha\Omega_{p,i},
\overline{\bar \pa_\gamma D_\alpha\Omega_{p,i}});\\
&(D_\alpha\Omega_{p,i},
\overline{\bar\pa_\gamma\pa_\gamma
D_\alpha\Omega_{p,i}})=-(\bar\pa_\gamma
D_\alpha\Omega_{p,i},
\overline{\bar \pa_\gamma D_\alpha\Omega_{p,i}}).
\end{align}
Inserting the above two equations into~\eqref{fourth},
we have
\begin{align}
\begin{split}
\frac{\pa^2 }{\pa t_\gamma\pa \bar
t_\gamma}h_{\alpha\bar\alpha}
&=\sum_p(\ii)^{q-p+2}(-(R_p)_{j\bar i\gamma\bar\gamma})
(D_\alpha\Omega_{p,i},\overline{D_\alpha\Omega_{p,j}})\\
&+\sum_p(\ii)^{q-p+2}(\pa_\gamma D_\alpha\Omega_{p,i},
\overline{\pa_\gamma D_\alpha\Omega_{p,i}})\\
&-\sum_p(\ii)^{q-p+2}(\bar\pa_\gamma
D_\alpha\Omega_{p,i},
\overline{\bar \pa_\gamma D_\alpha\Omega_{p,i}}).
\end{split}
\end{align}

By~\cite[page 33, Proposition 4]{Gr}, the curvature of
$(g_p)_{i\bar j}$ is
\begin{equation}\label{sixth}
(R_p)_{i\bar j\gamma\bar\gamma}
=(\ii)^{q-p}(D_\gamma\Omega_{p,i},
\overline{ D_\gamma\Omega_{p,j}})
-(\ii)^{q-p} (\bar\pa_\gamma\Omega_{p,i},\overline
{\bar\pa_\gamma\Omega_{p,j}}).
\end{equation}

Let
\begin{equation}\label{seven}
\pa_\gamma D_\alpha\Omega_{p,i}
=A_{p\gamma\alpha i}+B_{p\gamma\alpha i},
\end{equation}
where $A_{p\gamma\alpha i}\in H^{p-2,q+2}$ and $B_{p\gamma\alpha
i}\in H^{p-1,q+1}$. Then
\begin{align}\label{nineth}
\begin{split}
&\frac{\pa^2}
{\pa t_\gamma\pa\bar t_\gamma}h_{\alpha\bar\alpha}=
-\sum_p(\ii)^{q-p+2} (R_p)_{j\bar i\gamma\bar\gamma}
(D_\alpha\Omega_{p,i},\overline{D_\alpha
\Omega_{p,j}})\\
&-|A_{p\alpha\gamma i}|^2+|B_{p\alpha\gamma i}|^2
-\sum_p(\ii)^{q-p+2}(\bar\pa_\gamma
D_\alpha\Omega_{p,i},
\overline{\bar \pa_\gamma D_\alpha\Omega_{p,i}}).
\end{split}
\end{align}

By Lemma~\ref{lemmaA1} and Lemma~\ref{lemmaA2}, let
\[
D_\alpha\Omega_{p,i}=(A_\alpha^p)_{il}\Omega_{p+1,l},
\]
for matrices $A^p_\alpha=(A^p_\alpha)_{il}$. Then we have
\[
\bar\pa_\alpha\Omega_{p,i}=\overline{(A_\alpha^{p-1})_{li}}
\Omega_{p-1,l}.
\]
Thus from (\ref{sixth}),
(\ref{nineth}), in terms of the matrices $A^p_\alpha$, we have
\begin{align}
\begin{split}\label{A14}
&{\frac{\pa^2}{\pa t_\gamma\bar \pa
t_\gamma}}{h_{\alpha\bar\alpha}}
=\sum_{}(A_\gamma^p)_{jt}\overline{(
A^p_\gamma)_{it}}(A_\alpha^p)_{is}
\overline{(A_\alpha^p)_{js}}\\
&-\sum_{}
\overline{(A_\gamma^{p-1})_{tj}}{(A^{p-1}_\gamma)_{ti}}
(A_\alpha^p)_{is}\overline{(A_\alpha^p)_{js}}\\
&-\sum|\sum_t(A^p_\alpha)_{it}(A^{p+1}_\gamma)_{ts}|^2
+\sum|\sum_s(A_\alpha^p)_{is}
\overline{(A^{p}_\gamma)_{ts}}|^2\\
&+|B_{p\alpha\gamma i}|^2
\end{split}
\end{align}

The $H^{p-2,q+2}$ part $A_{p\gamma\alpha i}$ of
$\pa_\alpha D_\gamma\Omega_{p,i}$ is the same as  the
$H^{p-2,q+2}$ part of $\pa_\alpha\pa_\gamma\Omega_{p,i}$.
Thus
\[
A_{p\gamma\alpha i}=A_{p\alpha\gamma i}
\]
for $1\leq\alpha,\gamma\leq m$ and $p=0,\cdots,k$.
In terms of the matrices $A^p_\alpha$, we have
\[
A_{\gamma}^pA_\alpha^{p+1}=A_\alpha^pA_\gamma^{p+1}
\]
for $1\leq\alpha,\gamma\leq m$ and $p=0,\cdots,k$.
Using these commutative relations of the matrices
$A^p_\alpha$, from ~\eqref{A14}, we have

\begin{align}
\begin{split}\label{A15}
&\quad\quad \sum_{\alpha,\gamma}{\frac{\pa^2}{\pa t_\gamma\bar
\pa t_\gamma}}{h_{\alpha\bar\alpha}}=\sum_{p,\alpha\gamma,i}
|B_{p\alpha\gamma i}|^2\\&
+\sum_{p,\alpha,\gamma}{\rm Tr} \left[\left(
\overline{(A_\gamma^{p-1})^T}A_\alpha^{p-1}
-A_\alpha^p\overline{(A_\gamma^p)^T}\right) \overline{\left(
\overline{(A_\gamma^{p-1})^T}A_\alpha^{p-1}
-A_\alpha^p\overline{(A_\gamma^p)^T}\right)}^T\right].
\end{split}
\end{align}

Since the last term of the above expression is zero if and
only if $A^p_\alpha\equiv 0$ for any $\alpha$ and $p$,
there
exists an $\eps>0$, such that
\begin{equation}\label{A11}
{\frac {\pa^2}{\pa t_\gamma\pa\bar
t_\gamma}}{h_{\alpha\bar\alpha}}\geq \eps\sum_p{\rm
Tr}(A_\alpha^p\bar{(A_\alpha^p)}^T) \geq\eps
|h_{\alpha\bar\alpha}|^2.
\end{equation}

Finally, from~\eqref{impot},and~\eqref{A11}, we have
\[
\Delta f\geq \frac{\eps}{m} f^2-C_1 f.
\]
Since the Ricci curvature of $\omega_{WP}$ is lowerly
bounded and since $f$ is nonnegative,
the generalized maximum principal (cf.~\cite{CY1}) gives
\[
f\leq mC_1/\eps.
\]
and the theorem is thus proved.

\qed

%\bibliographystyle{abbrv}
%\bibliography{:::bib:bib}
%\bibliography{bib}

\end{document}